\documentclass{article}

\usepackage{epsfig}
\usepackage{rotating}
\def\gridbox#1/#2/#3{
\vbox to #1truecm{#3
\ifshowgrid\tickcount=0
  \loop\cgridw%
   \vbox to 0pt{\kern\tickcount truecm\hrule width#2truecm height\gridwidth\vss}
   \nointerlineskip \advance\tickcount by 1
   \ifdim\tickcount pt<#1pt\repeat 
  \hbox to 0pt{\tickcount=0\tick#1/\advance\tickcount by 1 %
 \loop\ifdim\tickcount pt<#2pt\nexttick1#1/\advance\tickcount by 1 \repeat\hss}
\else \vbox to 0pt{\hrule width#2truecm height0pt\vss}
\fi\vfil}\vfil}

\def\dpoint(#1,#2)#3{\vbox to 0pt{\kern#1
   \hbox{\kern#2{#3}}\vss}\nointerlineskip}

\def\point(#1,#2)#3{\dpoint(#1truecm,#2truecm){#3}}
\def\cpoint(#1,#2)#3{\setbox0=\hbox{#3}
   \dimen0=\ht0\advance\dimen0 by\dp0\divide\dimen0 by-2
   \advance\dimen0 by#1truecm\dimen1=\wd0\divide\dimen1 by-2
   \advance\dimen1 by#2truecm\dpoint(\dimen0,\dimen1){#3}}

\newcount\rmndr   
\def\rem#1#2{\rmndr=#1{}\divide\rmndr by#2{}%
\multiply\rmndr by-#2{}\advance\rmndr by #1}
\def\cgridw{\gridwidth\finegridw{}\rem\tickcount5{}
\ifnum\rmndr=0{}\gridwidth\roughgridw\fi}          

\def\tick#1/{\cgridw\vrule width\gridwidth height0pt depth#1truecm}
\def\nexttick#1#2/{\hbox to#1truecm{\hfil\tick#2/}}
\newcount\tickcount
\newdimen\finegridw\finegridw0.4pt\newdimen\roughgridw\roughgridw1.6pt
\newdimen\gridwidth
\newif\ifshowgrid \showgridtrue

\newtheorem{theorem}{Theorem}

\newtheorem{remark}{Remark}
\newtheorem{definition}{Definition}

\newtheorem{statement}{Statement}

\def\CC{{\mathcal C}}

\def\Re{{\mathcal R}e}
\def\Im{{\mathcal I}m}
 \title{An analogue of Gauss-Lucas theorem for a non convex sector of the complex plane}
\author{Bl. Sendov}
\date{}

\begin{document}

\maketitle
 
\begin{abstract} Let $S(\phi)= \{z:\;|\arg(z)|\geq \phi\}$ be a sector on the complex plane $\CC$. If $\phi\geq \pi/2$,
then $S(\phi)$ is a convex set and, according to the Gauss-Lucas theorem, if a polynomial $p(z)$ has all its zeros on $S(\phi)$, then the same is true for the zeros of all its derivatives. In this paper is proved that if the polynomial $p(z)$ is with real and non negative coefficients, then the same is true also for $\phi < \pi/2$, when the sector is not a convex set.


\end{abstract}

Let $P_n^+$ be the set of all algebraic polynomials of degree $n$ with real and non negative coefficients. Denote by 
$S(\varphi)=\{z:\; |\arg(z)|\geq\varphi\in [0,\pi],\;z\in\CC\}$
a sector on the complex plane $\CC$. Let $P_n^+(\varphi)$ be the set of all polynomials from 
$P_n^+$ with zeros on $S(\varphi),\; \varphi \in[0,\pi]$. 
 
\begin{theorem}[Sector theorem]\label{sectortheorem} If $p(z)\in P_n^+(\varphi)$ and $n\geq2$, then $p'(z)\in P_{n-1}^+(\varphi)$. 
\end{theorem}

It is clear that Theorem \ref{sectortheorem} is true for $\varphi\in [\pi/2,\pi]$, according to the Gauss-Lucas theorem, as in this case the sector  $S(\varphi)$ is a convex set. It is easy to verify the Theorem \ref{sectortheorem} for $ n=2,3$, so we may suppose that $n\geq 4$.

As the zeros of a polynomial are continuous functions of its coefficients and vice verse, we may
suppose that the polynomials in $P_n^+$ has only strictly positive coefficients and simple zeros.

\begin{statement}\label{lem}  Let $p(z)\in P_n^+$ has no zeros on the ray $L(\theta)=\{z:\; z=te^{i\theta};\; t\geq 0\}$. Denote by $\Delta(p;\theta)$ the net change of $\arg(p(z))$ when $z$ traverse the ray $L(\theta)$ from 0 to $\infty$.  
 Then
 $$\Delta(p;\theta))= n\theta -2m\pi;\; m\in \{0,1,2,\dots\},$$ 
 where $m$ is the number of the zeros of $p(z)$ inside the sector $s(\theta)= \{z:\;z=re^{i\psi};\; \psi 
 \in [0,\theta]\}$.
 \end{statement}
 
 The following is geometrically clear:
 \begin{statement}\label{also} If $u(t)$ and $v(t)$ are real polynomials, the zeros of which interlace, then the zeros of $U(t)=a u(t)+ b v(t)$ and $V(t)= c u(t)- d v(t)$, where $a,b,c,d >0$ also interlace.
 \end{statement}   
 \begin{definition}\label{weak} Let $p(x)$ and $q(x)$ be two real polynomials which degrees are equal or differ by 1. We say that the real zeros of these polynomials {\bf weakly interlace} if there is no two consecutive intervals defined by three consecutive zeros of the one polynomial, such that neither of them contained atleast one zero of the other polynomial and vice verse.
\end{definition}  

It is easy to see that:

\begin{statement}\label{weak1} Let $q_1(x), q_2(x), g_1x)$ and $g_2(x)$ be four real polynomials. If the zeros of the pairs $(q_1,q_2), (q_1,g_1)$ and $(q_2,g_2)$ interlace, then the zeros of the pair $(g_1,g_2)$ weakly interlace. 
\end{statement} 

\begin{remark}\label{remark} It is clear that in Statement \ref{also}, the word interlace may be replaced by weakly interlace.
\end{remark}
Consider the polynomials in $t$
$$g_1(t)= \Im\Bigl(p(te^{i\theta})\Bigr)=\sum_{k=1}^na_k t^k \sin k\theta,$$ 
$$g_2(t)= \Re\Bigl(p(te^{i\theta})\Bigr)=\sum_{k=0}^na_k t^k \cos k\theta,$$
$$g'_1(t)= \sum_{k=1}^n k a_k t^{k-1} \sin k\theta,\;\; g'_2(t)= \sum_{k=1}^n k a_k t^{k-1} \cos k\theta,$$ 
$$h_1(t)= \Im\Bigl(p'(te^{i\theta})\Bigr)=\sum_{k=2}^n k a_k t^{k-1} \sin (k-1)\theta,$$
$$h_2(t)= \Re\Bigl(p'(te^{i\theta})\Bigr)=\sum_{k=1}^n k a_k t^{k-1} \cos (k-1)\theta.$$

We have
$$g'_1(t)=h_1(t)\cos\theta + h_2(t)\sin\theta,\;\;\;\;g'_2(t)=h_2(t)\cos\theta -h_1(t)\sin\theta,$$
hence
\begin{equation}\label{equ1}
h_1(t)= g_1'(t)\cos\theta -g_2'(t)\sin\theta,\;\;\;\;h_2(t)=  g_1'(t)\sin\theta +g_2'(t)\cos\theta.
\end{equation}

{\bf Proof of Theorem \ref{sectortheorem}.} 
  Suppose that there exists a polynomial $p(z)$ with $p'(\zeta_0)=0$, where $\zeta_0= re^{i\varphi}$ be the zero of $p'(z)$ with the smallest positive argument $\varphi\in (0,\pi/2)$, which do not satisfy Theorem \ref{sectortheorem}. Then there exists $\theta>\varphi$, such that $p(z)\in P_n^+(\theta+\varepsilon);\; \varepsilon>0$ and $p'(z)$ has not zeros over $L(\theta)$. Without loss of the generality, we may suppose that $\zeta_0$ is the only zero of $p'(z)$ with positive argument less than $\theta$. 

From Statement \ref{weak} follows that 
\begin{equation}\label{theta}  
\Delta(p;\theta)= n\theta,\;\;\;\;\Delta(p';\theta)= (n-1)\theta - 2\pi.
\end{equation}   
 
 The number of the variations of the coefficients of $g_1(t)$ is equal to $[n\theta/\pi]$. By the Descartes' Rule of Signs, the number of the positive zeros of $g_1(t)$ is equal to the number of the variations of its coefficients or less with a even number. It is clear, that to change the  $\arg(p(z))$ by $\pi$, the value of $g_1(t)$ has to go from a zero to another one, as $g_1(0)=0$. Hence, from the first equation (\ref{theta}) follows that the number of the real positive zeros of $g_1(t)$ is exactly $[n\theta/\pi]$. This property is used also in \cite[p. 82]{RWKD}. 

 Let $0=t_0<t_1<t_2<\cdots <t_m;\; m=[n\theta/\pi]$ be the zeros of $g_1(t)$. To every interval $[t_{k-1},t_k];\;k=1,2,\dots,m$ correspond an increase of $\arg(p(z))$ with $\pi$. This is possible, if in every such interval, there exists even number of zeros of $g_2(t)$.
The  number of the variations of the coefficients of $g_2(z)$ is equal or with one bigger from this of $g_1(z)$, but $g_2(0)=a_0>0$, hence the number of the zeros of $g_2(z)$ is equal ore with one less from this of $g_1(z)$. From this follows that the  zeros of $g_1(t)$ and $g_2(t)$ interlace.   
 
From the Rolle's theorem and the Descartes' Rule of Signs follows that the non negative zeros of $g'_1(t)$ interlace withe the non negative zeros of $g_1(t)$ and the non negative zeros of $g_2'(t)$ interlace withe the non negative zeros of $g_2(t)$.
From Statement (\ref{weak1}) follows that the non negative zeros of the polynomials $g_1'(t)$ and $g_2'(t)$ weakly interlace.   Hence, from (\ref{equ1}) and Remark \ref{remark} follows
\begin{statement}\label{iterl} The non negative zeros of the polynomials $h_1(t)$ and $h_2(t)$ weakly interlace.  
\end{statement}

We already see that, if the non negative zeros of the polynomials $h_1(t)$ and $h_2(t)$ interlace, then $\Delta(p';\theta)= (n-1)\theta$, which contradicts to the second equation (\ref{theta}). Then, there is an interval $\Delta'$ between two consecutive  zeros of $h_1(t)$, which is free from the zeros of $h_2(t)$. There are two possibilities:

1) The interval $\Delta'$ is the last between two consecutive  zeros of $h_1(t)$. In this situation the net change of the $\arg(p'(z))$, when $t$ traverse $\Delta'$ is 0 and $\Delta(p',\theta)=(n-1)\theta-\pi$, which contradicts  the second equation (\ref{theta}).  
 
2) The interval $\Delta'$ is followed by the interval $\Delta''$, which contains a zero of $h_2(t)$. Then the net change of the $\arg(p'(z)$, when $t$ traverse $\Delta''$ is $-\pi$. We are loosing $\pi$ for the interval $\Delta'$ and $2\pi$ for the interval $\Delta''$. Hence $\Delta(p',\theta)=(n-1)\theta -3\pi$, which also contradicts  the second equation (\ref{theta}) and completes the proof.   

\begin{remark} The Sector theorem was already published \cite{BS} with a rather long proof. But the reason to publish a new short proof is mainly that the long one turned not to be correct. 
\end{remark}
  
{\large \bf Acknowledgments.} This work has been partly supported by the Bulgarian National Science Fund under project FNI I 02/20 "Efficient Parallel Algorithms for Large-Scale Computational Problems"

\bigskip 

\bigskip 

Address: Bl. Sendov

Acad. G. Bonchev str., Bl. 25

1113 Sofia, Bulgaria

E-mail: sendov2003@yahoo.com

\end{document}